

\font\seventeensans=cmss17
\font\smallcaps=cmcsc10
\font\twelvebf=cmbx12
\input amssym.def
\input amssym.tex

\def\mathbb{\Bbb}
\def\R{{\Bbb R}}
\def\N{{\Bbb N}}

\def\C{{\Bbb C}}

\def\Rn{{\R^n}}

\magnification1200

\null\vfill

\centerline{\seventeensans On cubic-linear polynomial mappings}

\bigskip\bigskip

\centerline{\smallcaps Gianluca Gorni}

\smallskip

\centerline{\it Universit\`a di Udine, Italy}

\bigskip

\centerline{\smallcaps Gaetano Zampieri}

\smallskip

\centerline{\it Universit\`a di Torino, Italy}

\bigskip\bigskip

\vfill

\noindent
{\bf Abstract.} In the field of the Jacobian conjecture it is
well-known after Dru\.zkowski that from a polynomial
``cubic-homogeneous'' mapping we can build a higher-dimensional
``cubic-linear'' mapping and the other way round, so that one of
them is invertible if and only if the other one~is. We make this
point clearer through the concept of ``pairing'' and apply it to
the related conjugability problem: one of the two maps is
conjugable if and only if the other one is; moreover, we find
simple formulas expressing the inverse or the conjugations of one
in terms of the inverse or conjugations of the other. Two
nontrivial examples of conjugable cubic-linear mappings are
provided as an application.

\vfill

\bigskip

\noindent
\smallcaps Gianluca Gorni,
{\it Universit\`a di Udine, Dipartimento di Matematica e
Informatica, via delle Scienze 208, 33100 Udine, Italy.}
{\tt Gorni@dimi.uniud.it}

\smallskip

\noindent
\smallcaps Gaetano Zampieri,
{\it Dipartimento di Matematica, via Carlo Alberto 10,
10123 Torino, Italy.} {\tt Zampieri@dm.unito.it}
\rm


\bigskip

\noindent
We are grateful to the Mathematics Department of the
Jagiellonian University, Krak\'ow, and in particular to Ludwik
Dru\.zkowski, for the invitation, the warm hospitality and the
stimulating discussions we had in a recent visit.

\smallskip

\noindent
We thank the University of Udine and the CNR-GNFM for providing us
computing facilities, that we put to heavy use in the early
stages of this research.

\eject

\pageno=2

\headline{\sevenrm On cubic-linear polynomial mappings\hfil
Gianluca Gorni and Gaetano Zampieri}

\def\Cn{\C^n}
\def\CN{\C^N}
\def\X{{\Bbb X}}

\def\det{\mathop{\rm det}}
\def\dom{\mathop{\rm dom}}
\def\ker{\mathop{\rm ker}}
\def\range{\mathop{\rm range}}

\def\mapright#1{\smash{\mathop{\longrightarrow}\limits^{#1}}}
\def\mapleft#1{\smash{\mathop{\longleftarrow}\limits^{#1}}}
\def\mapdown#1{\Big\downarrow
  \rlap{$\vcenter{\hbox{$\scriptstyle#1$}}$}}
\def\mapup#1{\Big\uparrow
  \rlap{$\vcenter{\hbox{$\scriptstyle#1$}}$}}
\def\bignearrow{\mathrel{\vcenter{\hbox{$\diagup$\kern-1.75pt
  \raise7.95pt\hbox{$\nearrow$}}}}}

\def\qed{
  \vbox{\hrule\hbox to7.8pt{\vrule height7pt
  \hss\vrule height7pt}\hrule}}

\centerline{\twelvebf 1. Introduction}

\bigskip\bigskip

The following conjecture was essentially originated by Keller~[14]
in~1939:

\bigskip

\noindent
{\bf Jacobian Conjecture}. \sl For all~$n\in\N$, if
$f\colon\Cn\to\Cn$ has polynomial components and the Jacobian
determinant $\det f'(x)$ is a nonzero constant throughout~$\Cn$,
then $f$~is a polynomial automorphism of~$\Cn$, that is, a
bijective polynomial map with polynomial inverse.
\rm

\bigskip

There is a huge literature on this topic and also some wrong
proofs were published. A~basic paper on the subject is~[2] by
Bass, Connell and Wright. The recent proceedings of 
conference~[10], and in particular its first paper, by the editor 
van den Essen, are a good update on this research field, rich in
questions of different nature.

In everything that follows $\R$ or $\Rn$ can be substituted for
$\C$ and~$\Cn$ with only trifling adjustments. Before proceeding it
is convenient to establish some notations first: if $x,y\in\Cn$ we
will write $x*y$ for the componentwise product of the two vectors:
$x*y:=(x_1y_1,x_2y_2,\ldots, x_ny_n)\allowbreak\in\Cn$. The powers
with respect to this multiplication will be denoted by
$x^{*2},x^{*3},\ldots$ The symbol $I_n$ will be identity mapping
(or matrix) in~$\Cn$. The minus signs in formulas~(1.1) and~(1.2)
below may seem odd but will later simplify some expressions in
Section~4.

\bigbreak

{\bf Definition~1.1.} \sl A mapping $f\colon\Cn\to\Cn$ will be
called ``cubic-homogeneous'' if there exists a trilinear symmetric
function $g\colon\Cn\times\Cn\times\Cn\to\Cn$ such that
$$f(x)=x-g(x,x,x)
  \qquad\hbox{for all }x\in\Cn\,,
  \eqno(1.1)$$
It will be called ``cubic-linear'' if there exists an $n\times n$
matrix~$A$ such that
$$f(x)=x-(Ax)^{*3}
  \qquad\hbox{for all }x\in\Cn\,.
  \eqno(1.2)$$
Cubic-homogeneous and cubic-linear mappings with constant Jacobian
determinant will be called ``Yagzhev maps'' and ``Dru\.zkowski
maps'' respectively.
\rm

\bigbreak

Two classical ``reduction'' results bear in particular on the
present paper. They restrict the class of polynomial functions over
which it is sufficient to concentrate the attention in order to
prove or disprove the full conjecture: the first reduction was to
Yagzhev maps (Yagzhev~[18] and independently
Bass-Connell-Wright~[2]) and the second was to the smaller class of
the Dru\.zkowski maps (Dru\.zkowski~[7]).

A side issue of the Jacobian conjecture was introduced in~[5]:
given a parameter $\lambda\in\C\setminus\{0,1\}$ and a polynomial
mapping~$f\colon\Cn\to\Cn$ such that $f(0)=0$, $f'(0)=I_n$, the
problem is to find a global analytic {\it conjugation}, i.e., an
invertible analytic function $k_\lambda\colon\Cn\to\Cn$ such that
the following diagram commutes:
$$
  \matrix{\C^n&\mapleft{k_\lambda}&\Cn\cr
          \mapdown{\lambda F}&&\mapdown{\lambda I_n}\cr
          \C^n&\mapleft{k_\lambda}&\Cn\cr}
  \eqno(1.3)$$
and with again the ``normalizing'' conditions $k_\lambda(0)=0$,
$k'_\lambda(0)=I_n$. We will also be handling functions
$k_\lambda$ that are like conjugation, except for either being
defined only on a neighbourhood of~$0\in\Cn$ or for possibly
failing to be invertible: these will be called {\it
pre-conjugations}.

\goodbreak

If a mapping~$f\colon\Cn\to\Cn$ admits a conjugation for some
$\lambda$, then the discrete dynamical system of the backward and
forward iterates of~$\lambda f$ is ``trivial''. In particular, $f$
is invertible, and this observation was in fact the original
motivation for raising the problem. There was some hope to prove
that the Jacobian conjecture was true by proving first that all
polynomial maps in a suitable class are conjugable.

With the work that has been done in the past few years on
conjugations, and specially after the examples found by A.~van
den~Essen and E.~Hubbers, the hope to possibly prove the Jacobian
conjecture through conjugations has dimmed, although we cannot rule
it out yet. What is still very well possible is that a
counterexample to the Jacobian conjecture may be found as a
by-product of research in conjugability. One unexpected and
encouraging by-product, in a seemingly unrelated area, is already
here: it is a very simple and elegant counterexample to
Markus-Yamabe conjecture in dimension~$\ge3$, due to  A.~Cima, 
A.~van den~Essen, A.~Gasull, E.~Hubbers and  F.~Ma\~nosas~[4].

The state of the art in the conjugation business is as follows.
There are normalized polynomial automorphisms that are not
conjugable: the earliest one was given in~[11], and it is of
``quintic-homogeneous'' form, but we have later realized that also
the old and well-known map in two dimensions
$$f\colon{x\choose y}\mapsto{x+(y+x^2)^2\choose y+x^2}
  \eqno(1.4)$$
is a counterexample (three fixed points for $\lambda f$ outside the
origin are quickly found when $\lambda\in\C\setminus\{0,1\}$). We
do not know yet if all Yagzhev maps are necessarily conjugable,
although a theorem in~[12] seems to make it hard to find
counterexamples. Yagzhev maps have been given for which global
conjugations exist that are analytic but not polynomial
themselves~([10, page~231] and~[12]). These last examples have in
turn taught us a lesson on the Jacobian conjecture, namely, on the
structure of the local inverse of Yagzhev maps (see~[13]).

The present paper was born out of the effort to find whether a
particular Dru\.zkowski mapping (Example~6.1 below) was conjugable.
We have discovered that there is a strong link between the
conjugability of a Dru\.zkowski map, or, more generally, of a
cubic-linear map~$F$, and the conjugability of a certain
lower-dimensional cubic-homogeneous map~$f$. The exact relation
between $F$ and~$f$ is described as follows:

\bigbreak

{\bf Definition~1.2.} \sl Given a cubic-homogeneous mapping
$f\colon\Cn\to\Cn$ and a cubic-linear mapping $F\colon
\CN\to\CN$, $F(X):=X-(AX)^{*3}$, with~$N>n$, we will say that $f$
and~$F$ are ``paired'' through the matrices $B$ and~$C$
(of~dimensions $n\times N$ and $N\times n$ respectively) if $\ker
A=\ker B$ and the following diagrams commute:
$$
  \matrix{\CN&\mapright B&\Cn\cr
  \mapup C&\bignearrow_{\rlap{{$\scriptstyle I_n$}}}&\cr
  \Cn&&\cr}
  \qquad\qquad
  \matrix{\CN&\mapleft C&\Cn\cr
  \mapdown F&&\mapdown f\cr
  \CN&\mapright B&\Cn\cr}
  \eqno(1.5)$$
that is, $BC=I_n$ and $f(x)=BF(Cx)$ for all~$x\in\Cn$. 
\rm

\bigbreak

We will prove that when two maps are paired, one has a conjugation
if and only if the other has. But the symmetry extends to
invertibility too: in fact it turns out that the pairing concept
underlies and somewhat elucidates Dru\.zkowski's reduction theorem,
in particular the way it is proved in~[9], and~[10, page~11].
The main results of this paper can be summed up in the following
theorem, which shows that the problems of both invertibility and
conjugation have the same answers if two maps are paired.

\bigbreak

{\bf Theorem~1.3.} \sl Every cubic-homogeneous map can be paired to
a cubic-linear map and vice versa. Moreover, if $f$ and~$F$ are
paired, each of the following properties for one of the two
mappings implies the same property for the other, for a
given~$\lambda\in\C\setminus\{0\}$, $|\lambda|\ne1$:

\item{1.}one-to-one,
\item{2.}onto,
\item{3.}invertible with polynomial inverse,
\item{4.}constant Jacobian determinant,
\item{5.}existence of a global pre-conjugation,
\item{6.}the global pre-conjugation is onto,
\item{7.}the global pre-conjugation is one-to-one,
\item{8.}the global pre-conjugation is a polynomial map,
\item{9.}the conjugation has a polynomial inverse.

\noindent
The existence of global pre-conjugations is guaranteed when
$|\lambda|>1$, and also when $f$ and~$F$ are invertible.

\rm

\bigbreak

Here is an assortment of formulas connecting $F,f$ and
the respective (globally defined) pre-conjugations $K_\lambda,
k_\lambda$, that are true for all $x,y\in\Cn$, $X,Y\in\CN$,
whenever every single piece just makes sense:
$$\eqalign{&f(BX)=BF(X)\,,\qquad
  \det f'(x)=\det F'(Cx)\,,\qquad
  \det F'(X)=\det f'(BX)\,,\cr
  &f^{-1}(y)=BF^{-1}(Cy)\,,\qquad
  F^{-1}(Y)=Y-F\bigl(Cf^{-1}(BY)\bigr)+Cf^{-1}(BY)\,,\cr
  &k_\lambda(x)=BK_\lambda(Cx)\,,\qquad
  k_\lambda(BX)=BK_\lambda(X)\,,\cr
  &k^{-1}_\lambda(y)=BK_\lambda^{-1}(Cy)\,,\qquad
  K_\lambda^{-1}(Y)=Y-K_\lambda\bigl(Ck_\lambda^{-1}(BY)\bigr)+
  Ck_\lambda^{-1}(BY)\,.\cr}
  \eqno(1.6)$$
The first two rows (and in particular the formula for $F^{-1}$ in
terms of~$f^{-1}$) were somehow implicit in the treatment of~[9]
and~[10, page~11], but they were hidden beneath layers of changes of
variables.
  
The rest of the paper is organized as follows. Section~2 concerns
the existence of pairing and some basic properties. Section~3 is
about invertibility. Section~4 is an introduction to
pre-conjugations, in the simpler cubic-homogeneous setting and with
a much easier proof than in~[5], and not relying on Poincar\'e's
theorem~[1, Sec.~25] either. Section~5 shows how pairing behaves
under conjugation. Section~6 illustrates two examples: the first is
Dru\.zkowski's example~7.8 from~[9] in dimension~15, which turns
out to be conjugable through a polynomial automorphism; in the end
we compute a Dru\.zkowski pairing to van den~Essen's example
from~[10, page~231], thus producing a new Dru\.zkowski map in
dimension~16 for which global analytic conjugations exist which are
not polynomial.

\vfill\eject


\centerline{\twelvebf 2. Pairing a cubic-homogeneous mapping}
\smallskip
\centerline{\twelvebf    to a cubic-linear one, and vice versa} 

\bigskip

{\bf Proposition~2.1.} \sl Let $f\colon\Cn\to\Cn$ be a
cubic-homogeneous mapping. Then there exist $N>n$
and linear maps $A\colon\CN\to\CN$, $B\colon\CN\to\Cn$,
$C\colon\Cn\to\CN$ such that $f$ is paired to the cubic-linear
mapping $F(X):=X-(AX)^{*3}$ through~$B$ and~$C$.
\rm

\medbreak

{\bf Proof.} (To follow the steps of this proof it may
help to look at the last example of Section~6, where they are
carried out in some detail on a nontrivial mapping~$f$). Thanks to
the algebraic identities (see~[9])
$$\eqalign{
  ab^2={}&{(a+b)^3+(a-b)^3-2a^3\over6}\,,\cr
  abc={}&{(a+b+c)^3+(a-b-c)^3-
  (a+b-c)^3-(a-b+c)^3\over24}\,,\cr}
  \eqno(2.1)$$
we can write every third-degree monomial appearing in the
components of~$f$ as a linear combination of cubic powers of linear
forms of~$x$. Build a matrix~$D_0$ by piling up in some order all
the 1-row matrices representing these linear forms. Do not forget
to insert the projections corresponding to the monomials such
as~$x_1^3$, that are cubic powers from the start. Next, build the
matrix $B_0$ that combines the cubes of those linear forms so that
$$f(x)=x-B_0(D_0x)^{*3}
  \qquad\hbox{for all }x\in\Cn\,.
  \eqno(2.2)$$
The matrix~$B_0$ has the same dimensions as the transpose
of~$D_0$. By adding null columns to~$B_0$ and an equal number of
null rows to~$D$ we can assume that the number of columns of~$B_0$
is~$>n$. The matrix $B_0$ may not yet be the~$B$ of the
statement, because it need not be of full rank. But this problem
is easily remedied by adding a few columns to~$B_0$ and the same
numbers of null rows to~$D_0$. For example we can add a
$n\times n$ identity matrix at the right end of~$B_0$ and a
$n\times n$ null matrix to the bottom of~$D_0$. In a similar
manner we can arrange that $D_0$ has full rank too. Call $B,D$ the
resulting matrices and $N$ the number of columns of~$B$. We have
that $N>n$ and
$$f(x)=x-B(Dx)^{*3}
  \qquad\hbox{for all }x\in\Cn\,.
  \eqno(2.3)$$
Let $C$ be any right-inverse of~$B$, i.e., an~$N\times n$
matrix such that $BC=I_n$. What we are still missing is an
$N\times N$ matrix~$A$ that shares the same kernel as~$B$ and such
that $f(x)=BF(Cx)$, where $F$ is defined as $F(X):= X-(AX)^{*3}$.
Let~$M$ be a matrix whose columns form a basis of the kernel
of~$B$. If~we are content for the time being to relax the equality
of the kernels into the inclusion $\ker A\supset\ker B$, this weaker
condition in terms of~$M$ translates as $AM=0$. On the other hand,
if we impose that $AC=D$, we will be able to write $f(x)=
x-B(ACx)^{*3}=B(Cx-(ACx)^{*3})=BF(Cx)$. The two equations $AM=0$
and $AC=D$ can be combined as
$$A(C\mid M)=(D\mid0)\,,
  \quad\hbox{which solves for $A$ as}\quad
  A=(D\mid0)(C\mid M)^{-1}\,,
  \eqno(2.4)$$
where $(C\mid M)$ is the matrix formed by joining the two blocks of
columns of~$C$ and of~$M$, and $(D\mid0)$ similarly. The matrix
$(C\mid M)$ is indeed invertible because the range of~$C$ is a
complement to the kernel of~$B$, since $BC=I_n$. The proof is
complete if we notice that with this choice of~$A$ the kernel
of~$B$ is equal to, and not merely contained in, the kernel of~$A$,
because the rank of~$A$ is the same as the rank of~$D$.
\qed

\bigbreak

The reverse procedure from a cubic-linear to a paired
cubic-homogeneous mapping is much easier. Throughout the rest of
this paper $A$ will be a fixed linear mapping $A\colon\CN\to\CN$,
that we will as usual identify with the matrix that represents it
with respect to the canonical basis of~$\CN$, and $F(X):=
X-(AX)^{*3}$ for~$X\in\CN$. The matrix $A$ will be assumed to be
singular, both because this is the case when the Jacobian
determinant is constant (that is, if we are dealing with what we
called Dru\.zkowski maps; see~[7]), and because the following
theory trivializes anyway when $A$ is invertible. Before
proceeding, take note of the following fact, that we will be using
over and over again.

\bigbreak

{\bf Proposition~2.2.} \sl If~$X\in\CN$ and $X_0\in\ker A$, then
$F(X+X_0)=F(X)+X_0$. In particular the differential satisfies
$F'(X)X_0=X_0$ and $F'(X+X_0)=F'(X)$.
\rm

\medbreak

{\bf Proof.} Obvious: $F(X+X_0)=X+X_0-(AX+AX_0)^{*3}=F(X)+X_0$.
\qed

\bigbreak

Let $n$ be the rank of~$A$ and $B\colon\CN\to\Cn$ be a linear
mapping with the same kernel as~$A$. In particular $B$ has full
rank, coinciding with the rank of~$A$. Let $C\colon\Cn\to\CN$ be a
right-inverse of~$B$, that is, a linear mapping such that
$BC=I_n$. A mapping~$f$ that is paired to $F$ through
$B$ and~$C$ is trivial to define:
$$f(x):=BF(Cx)=x-B\bigl(ACx\bigr)^{*3}
  \quad\hbox{for }x\in\Cn\,.
  \eqno(2.5)$$
A property of~$B$ and~$C$ that we will also be using all the time
without explicit reference is that
$$CBX-X\in\ker A=\ker B \qquad
  \hbox{for all }X\in\CN\,.
  \eqno(2.6)$$
The formula is true, because $B(CBX-X)=(BC)BX-BX=BX-BX=0$.

\bigbreak

{\bf Proposition~2.3.} \sl Once $B$ is given, the paired
mapping~$f$ defined in~(2.5) is independent of the choice of the
right-inverse~$C$, and it makes the following diagram commute:
$$
  \matrix{\CN&\mapright B&\Cn\cr
  \mapdown F&&\mapdown f\cr
  \CN&\mapright B&\Cn\cr}
  \eqno(2.7)$$
\rm

\medbreak

{\bf Proof.} Let $C,\tilde C$ be two right inverses of~$B$. Then
$Cx-\tilde Cx\in\ker A=\ker B$ for all~$x\in\Cn$, because
$B(Cx-\tilde Cx)= BCx-B\tilde Cx=x-x=0$. The paired mapping~$f$
does not depend on the choice of~$C$ because
$$BF(\tilde Cx)=B\Bigl(F(\tilde Cx)+
  \underbrace{Cx-\tilde Cx}_{\hbox to0pt{\hss
  $\scriptstyle\in\ker B=\ker A$\hss}}\Bigr)=
  BF\bigl(\tilde Cx+Cx-\tilde Cx\bigr)=BF(Cx)\,.
  \eqno(2.8)$$
As for diagram~(2.7), noticing that $ACB=A$,
$$\eqalign{
  f(BX)={}&BF(CBX)=B\Bigl(F(CBX)+
  \underbrace{X-CBX}_{\hbox to0pt{\hss
  $\scriptstyle\in\ker B=\ker A$\hss}}\Bigr)=\cr
  ={}&BF\bigl(CBX+X-CBX\bigr)=BF(X)\,.\cr}
  \eqno(2.9)$$
\qed

\bigbreak

{\bf Proposition~2.4.} \sl  For all~$x\in\Cn$, $X\in\CN$
we have $\det f'(x)=\det F'(Cx)$ and $\det F'(X)=\det f'(BX)$. In
particular $f$ has constant Jacobian determinant if and only if
$F$~has.
\rm

\medbreak

{\bf Proof.}  To study the Jacobian determinants of~$F$ and~$f$ it
is convenient to decompose first $\CN=(\range C)\oplus(\ker A)$
and to choose a basis of~$\CN$ whose first $n$ vectors are
the image through~$C$ of the canonical basis of~$\Cn$
(forming in particular a basis of the range of~$C$) and the
remaining ones are a basis of~$\ker A$. If on~$\Cn$ we
keep the canonical basis, the matrices representing $F'(X),C,B$
take the following forms, thanks also to Proposition~2.2,
$$F'(X)=\left(\matrix{
  R(X)&\vrule height15pt depth8pt&S(X)\cr
  \noalign{\hrule}
  0&\vrule height18pt depth8pt&I_n\cr}
  \right)\,,\qquad
  C=\left(\matrix{I_n\vrule height15pt
  depth8pt width0pt\cr\noalign{\hrule}0
  \vrule height15pt depth8pt width0pt\cr}
  \right)\,,\qquad
  B=\left(\matrix{I_n&\vrule height15pt depth8pt&0\cr}
  \right)\,,
  \eqno(2.10)$$
for matrices $R(X),S(X)$ of suitable dimensions. Now
$$\eqalign{
  \det f'(x)={}&\det BF'(Cx)C=\cr
  ={}&\det\left(\matrix{I_n&\vrule height15pt
  depth8pt&0\cr}
  \right)
  \left(\matrix{
  R(Cx)&\vrule height15pt depth8pt&N(Cx)\cr
  \noalign{\hrule}
  0&\vrule height18pt depth8pt&I_n\cr}
  \right)
  \left(\matrix{I_n\vrule height15pt
  depth8pt width0pt\cr\noalign{\hrule}0
  \vrule height15pt depth8pt width0pt\cr}
  \right)=\cr
  ={}&\det R(Cx)=\det F'(Cx)\,.\cr}
  \eqno(2.11)$$
Conversely,
$$\det F'(X)=\det F'\bigl(CBX+\underbrace{X-CBX}_{\hbox to0pt{\hss
  $\scriptstyle\in\ker B=\ker A$\hss}}\bigr)=
  \det F'(CBX)=\det f'(BX)\,.
  \eqno(2.12)$$
\qed

\vfill\eject


\centerline{\twelvebf 3. Inverses of paired mappings} 

\bigskip
\bigbreak

{\bf Proposition~3.1.} \sl If $F$ is one-to-one, so is~$f$.
If~$F$ is onto, so is~$f$. If~$F$ is a bijection, then so is~$f$,
and $f^{-1}(y)=BF^{-1}(Cy)$ for all~$y\in\Cn$, that is,
the following diagram commutes:
$$
  \matrix{\CN&\mapright B&\Cn\cr
  \mapup{F^{-1}}&&\mapup{f^{-1}}\cr
  \CN&\mapleft C&\Cn\cr}
  \eqno(3.1)$$
In particular, if $F^{-1}$ is a polynomial mapping, so
is~$f^{-1}$, and the degree of~$f^{-1}$ is not higher than the
degree of~$F^{-1}$.
\rm

\medbreak

{\bf Proof.} Suppose that $F$ is one-to-one. Let $x_0,x_1\in\Cn$.
Then
$$\eqalign{
  f(x_1)=f(x_2)\qquad\Longrightarrow{}&\qquad
  BF(Cx_1)=BF(Cx_2)\cr
  \Longrightarrow{}&\qquad
  F(Cx_1)-F(Cx_2)=:X_0\in\ker B=\ker A\cr
  \Longrightarrow{}&\qquad
  F(Cx_1-X_0)=F(Cx_2)\cr
  \Longrightarrow{}&\qquad
  Cx_1-X_0=Cx_2\cr
  \Longrightarrow{}&\qquad
  \range C\ni C(x_1-x_2)=X_0\in\ker A=\ker B\cr
  \Longrightarrow{}&\qquad
  C(x_1-x_2)=X_0=0\cr
  \Longrightarrow{}&\qquad
  x_1=x_2\,.\cr}
  \eqno(3.2)$$
Suppose that $F$ is onto. For a given $y\in\Cn$, we have to
prove that $y$ is in the range of~$f$. Let
$X\in\CN$ be such that $F(X)=Cy$. Then
$$\eqalign{
  \range f\ni{}&
  f(BX)=BF\bigl(CBX\bigr)=
  B\Bigl(F\bigl(CBX\bigr)+
  \underbrace{X-CBX}_{\hbox to0pt{\hss
  $\scriptstyle\in\ker B=\ker A$\hss}}\Bigr)=\cr
  ={}&B\Bigl(F\bigl(CBX+X-CBX\bigr)\Bigr)=
  BF(X)=BCy=y\,.\cr}
  \eqno(3.3)$$
Finally, when $F$ is a bijection, the vector~$X$ in~(3.3) is simply
$F^{-1}(Cy)$, which proves the first formula for the inverse.
\qed

\bigbreak

{\bf Proposition~3.2.} \sl If $f$ is one-to-one, so is~$F$.
If~$f$ is onto, so is~$F$. If~$f$ is a bijection, then so
is~$F$, and for all $Y\in\CN$
$$F^{-1}(Y)=Y+\Bigl(ACf^{-1}(BY)\Bigr)^{*3}=
  Y-F\bigl(Cf^{-1}(BY)\bigr)+Cf^{-1}(BY)
  \eqno(3.4)$$
In particular, if~$f^{-1}$ is a polynomial mapping, so
is~$F^{-1}$, and the degree of~$F^{-1}$ is at most three times the
degree of~$f^{-1}$.
\rm

\medbreak

{\bf Proof.} Suppose that $f$ is one-to-one and let
$X_1,X_2\in\CN$. Then
$$\eqalign{
  F(X_1)={}&F(X_2)\quad\Longrightarrow\cr
  \Longrightarrow{}&\quad
  F\Bigl(CBX_1+
  \underbrace{X_1-CBX_1}_{\hbox to0pt{\hss
  $\scriptstyle\in\ker B=\ker A$\hss}}\Bigr)=
  F\Bigl(CBX_2+
  \underbrace{X_2-CBX_2}_{\hbox to0pt{\hss
  $\scriptstyle\in\ker B=\ker A$\hss}}\Bigr)\cr
  \Longrightarrow{}&\quad
  F(CBX_1)+X_1-CBX_1=F(CBX_2)+X_2-CBX_2 \quad\hbox{(*)}\cr
  \Longrightarrow{}&\quad
  BF(CBX_1)=BF(CBX_2)\cr
  \Longrightarrow{}&\quad
  f(BX_1)=f(BX_2)\cr
  \Longrightarrow{}&\quad
  BX_1=BX_2\qquad\hbox{(using formula * above)}\cr
  \Longrightarrow{}&\quad
  X_1=X_2\,.\cr}
  \eqno(3.5)$$
Suppose that $f$ is onto and let~$Y\in\CN$. Let~$x\in\Cn$ be such
that $f(x)=BY$. Then
$$\eqalign{
  \range F\ni{}&
  F\bigl(Y+(ACx)^{*3}\bigr)=\cr
  ={}&
  F\Bigl(\underbrace{Y-CBY}_{\hbox to0pt{\hss
  $\scriptstyle\in\ker A$\hss}}+
  \underbrace{CBY}_{\hbox to0pt{\hss
  $\scriptstyle=Cf(x)$\hss}}+(ACx)^{*3}\Bigr)=\cr
  ={}&
  Y-Cf(x)+F\Bigl(
  \underbrace{Cf(x)-F(Cx)}_{\hbox to0pt{\hss
  $\scriptstyle\in\ker B=\ker A$\hss}}+
  \underbrace{F(Cx)+(ACx)^{*3}}_{\hbox to0pt{\hss
  $\scriptstyle=Cx$\hss}}\Bigr)=\cr
  ={}&
  Y-Cf(x)+Cf(x)-F(Cx)+F(Cx)=\cr
  ={}&Y\,.\cr}
  \eqno(3.6)$$
Assume finally that $f$ is a bijection. Then we can write
$x=f^{-1}(BX)$ in~(3.6) and get the first formula for the inverse.
The second expression is a simple consequence:
$$\eqalign{
  F^{-1}(Y)={}&
  Y+\Bigl(ACf^{-1}(BY)\Bigr)^{*3}=\cr
  ={}&Y-\biggl(Cf^{-1}(BY)
  -\Bigl(ACf^{-1}(BY)\Bigr)^{*3}\biggr)+Cf^{-1}(BY)=\cr
  ={}&Y-F\bigl(Cf^{-1}(BY)\bigr)+Cf^{-1}(BY)\,.\cr}
  \eqno(3.7)$$
\qed

\vfill\eject


\centerline{\twelvebf 4. Pre-conjugations for
                         cubic-homogeneous mappings} 

\bigskip

{\bf Proposition~4.1.} \sl Let $\X$ be a complex Banach space,
$\gamma\colon\X\times\X\times\X\to\X$ be a continuous trilinear
symmetric form, and define the function $\varphi\colon\X\to\X$ as
$\varphi(x):=x-\gamma(x,x,x)$. Then for any $\lambda\in\C
\setminus\{0\}$, with $|\lambda|\ne1$, there exists an analytic
function $\kappa_\lambda$ defined in a neighbourhood of $0\in\X$
and with values in~$\X$, such that
$$\eqalign{
  &\kappa_\lambda(0)=0\,,\quad \kappa_\lambda'(0)=I_\X
  \quad\hbox{(the identity operator on $\X$) and}\cr
  &\lambda\varphi(\kappa_\lambda(y))=\kappa_\lambda(\lambda y)
  \quad\hbox{ for all $y\in\X$ such that }
  y,\lambda y\in\dom \kappa_\lambda\,.\cr}
  \eqno(4.1)$$
The function $\kappa_\lambda$ is unique, in the sense that any two
functions with the same property must agree in a neighbourhood of
the origin. If we denote by $\Psi_m$ the homogeneous term of
degree~$m$ in the Taylor series $\kappa_\lambda=\sum_{m\ge0}\Psi_m$
of~$\kappa_\lambda$ centered in the origin (ignoring the dependence
on~$\lambda$), the following recursive formulas hold:
$$\eqalign{
  &\Psi_0(y):=0\,,\qquad\Psi_1(y):=y\,,\cr
  &\Psi_m={1\over1-\lambda^{m-1}}
  \sum_{p+q+r=m \atop 0\le p,q,r<m}
  \gamma\bigl(\Psi_p,\,\Psi_q,\,\Psi_r\bigr)\,,
  \quad\hbox{for }m\ge2\,.\cr}
  \eqno(4.2)$$
If either $|\lambda|>1$ or $\varphi$ is invertible, then the
function $\kappa_\lambda$ is defined and analytic on the whole
of~$\X$. Finally, if $\X$ is finite-dimensional and the Jacobian
determinant of~$\varphi$ is constant, then the same happens
to~$\kappa_\lambda$ on any connected open neighbourhood of the
origin (both constants must be~1, of course, because
$\varphi'(0)=\kappa_\lambda'(0)=I_\X$).
\rm

\medbreak

{\bf Proof.}
Uniqueness of $\kappa_\lambda$ and the recursive relations~(4.2)
are obtained as in~[13] simply by substitution of~$\kappa_\lambda=
\sum_m\Psi_m$ into the conjugation formula $\lambda
\varphi(\kappa_\lambda(y))=\kappa_\lambda(\lambda y)$, using the
multilinearity of~$\gamma$ and the homogeneity of~$\Psi_k$:
$$\lambda\sum_{m\ge0}\Psi_k-
  \lambda\sum_{p,q,r\ge0}\gamma\bigl(\Psi_p,\Psi_q,\Psi_r\bigr)=
  \sum_{m\ge0}\lambda^m\Psi_m\,,
  \eqno(4.3)$$
and then by grouping together the terms which are homogeneous of
the same degree. The initial conditions on $\Psi_0,\Psi_1$
cannot be derived from the conjugation relation, and are simply the
transcriptions of the normalizing conditions on~$\kappa_\lambda(0),
\kappa_\lambda'(0)$. The summation in~(4.2) can be restricted to
the $p,q,r$ strictly less than~$m$ because $\Psi_0=0$.
Observe that $\Psi_m=0$ when $m$~is even, a fact that we have
chosen not to highlight here, but that speeds up computations
sometimes.

\noindent
We have to prove that the series $\sum_k\Psi_k(y)$ converges when
$\|y\|$ is small enough. Write
$$\eqalign{
  a_m:={}&\sup_{\|y\|\le1}\bigl\|\Psi_m(y)\bigr\|\,,\cr
  \|\gamma\|:={}&\sup\Bigl\{\bigl\|\gamma(x,y,z)\bigr\|\;:\;
  \|x\|\le1,\,\|y\|\le1,\,\|z\|\le1\Bigr\}\,.\cr}
  \eqno(4.4)$$
The series $\sum\Psi_m(y)$ will converge whenever $\sum
a_m\|y\|^m<+\infty$. The following inequalities hold:
$$\eqalign{
  a_0={}&0\,,\qquad a_1=1\,,\cr
  a_m\le{}&{\|\gamma\|\over|1-\lambda^{m-1}|}
  \sum_{p+q+r=m}a_pa_qa_r\le
  {\|\gamma\|\over\bigl|1-|\lambda|\bigr|}
  \sum_{p+q+r=m}a_pa_qa_r\,.\cr}
  \eqno(4.5)$$
Then we see that $0\le a_m\le b_m$ for all~$m$, where $b_m$ is the
sequence defined by recursion as
$$b_0:=0\,,\quad b_1:=1\,,\quad
  b_m:=\alpha\sum_{p+q+r=m}b_pb_qb_r
  \quad\hbox{ for }m\ge2\,,\quad\hbox{where }
  \alpha:={\|\gamma\|\over\bigl|1-|\lambda|\bigr|}\,.
  \eqno(4.6)$$
If we define the one-variable (formal) power series $\mu(t):=
\sum b_mt^m$, we see that the function $\mu$ should verify the
relation
$$\mu(0)=0\,,\quad\mu'(0)=1\,,\qquad
  \mu(t)-\alpha\mu(t)^3=t
  \eqno(4.7)$$
for all~$t$ where $\mu(t)$ exists. This means that $\mu$
must be a local inverse of the complex variable function
$u\mapsto u-\alpha u^3$, around the origin, mapping 0 to~0. But we
very well know that such a local inverse exists and it is a power
series with a positive radius~$R$ of convergence. We could
estimate~$R$, if we wish, using Cardano's formula for cubic
equations. We conclude that the power series $\sum b_mt^m$ has
positive radius~$R$ of convergence. If~$\|y\|<R$ we have that
$\sum\|\Psi_m(y)\|\le\sum a_m\|y\|^m\le\sum b_m\|y\|^m=
\mu(\|y\|)<+\infty$. The local existence of~$\kappa_\lambda$ is
established.

\noindent
The fact that $\kappa_\lambda$ exists on the whole of~$\X$
if~$|\lambda|>1$ follows from the same simple argument
used in~[5]: the conjugation relation $\lambda\varphi
(\kappa_\lambda(y))= \kappa_\lambda(\lambda y)$ allows
us to define $\kappa_\lambda(\lambda y)$ whenever we know
$\kappa_\lambda(y)$, and the extensions that we obtain this way are
analytical.

\noindent
Similarly, when $\varphi$ is invertible, the conjugation relation can
be rewritten as $\kappa_\lambda(y)=\varphi^{-1}
(\kappa_\lambda(\lambda y)/\lambda)$, which allows us to extend
analytically the definition of~$\kappa_\lambda$ to the whole space
if $0<|\lambda|<1$.

\noindent
The derivative of the conjugation identity $\lambda
\varphi(\kappa_\lambda(y))= \kappa_\lambda(\lambda y)$ with
respect to~$y$ is $\lambda\varphi'(\kappa_\lambda(y))\allowbreak
\lambda\kappa_\lambda'(y)=\kappa_\lambda'(\lambda y)$. If~$\X=\Cn$
and $\varphi$ has constant Jacobian determinant, then this constant
is~1 because $\varphi'(0)=I_n$, and we deduce that
$\det\kappa_\lambda' (\lambda y)=\det \kappa_\lambda'(y)$.
If~$y\in\Cn\setminus\{0\}$ is close enough to the origin then
$\lambda^r y\in\dom \kappa_\lambda$ either for all $r\ge0$ or for
all $r\le0$, depending on whether $|\lambda|>1$ or $|\lambda|<1$.
In either case $\det \kappa_\lambda'$ has the same value along a
sequence of points containing~$y$ and with the origin as a cluster
point. Then $\det\kappa_\lambda'(y)=\det k'(0)=1$ because
$\kappa_\lambda'$ is continuous.
\qed

\bigbreak

{\bf Remark~4.2.} If we consider the local inverse of~$\varphi$
around the origin, the terms of its Taylor expansion
$\varphi^{-1}=\sum_m\Phi_m$ satisfy the same recursive relations
as the~$\Psi_m$, only with~$\lambda=0$ (see~[8]). It follows from
this with simple calculations that the~$\Psi_m$ are scalar
multiples of the corresponding $\Phi_m$ up to degree~5:
$$\Psi_1=\Phi_1\,,\qquad
  \Psi_3={1\over1-\lambda^2}\Phi_3\,,\qquad
  \Psi_5={1\over(1-\lambda^2)(1-\lambda^4)}\Phi_5\,.
  \eqno(4.8)$$
However the property fails from degree~7 onward. For example
$$\Psi_7(x)={1\over(1-\lambda^2)(1-\lambda^4)(1-\lambda^6)}
  \biggl(\Phi_7(x)+
  3\lambda^2
  \gamma\Bigl(\gamma(x,x,x),\gamma(x,x,x),x\Bigr)\biggr)\,.
  \eqno(4.9)$$

\vfill\eject


\centerline{\twelvebf 5. Conjugations of paired mappings} 

\bigbreak

In this section we will use the letters~$F,f,A,B,C$ with the same
meaning as in Section~2. The function $F$ can be expressed
as $F(X)=X-G(X,X,X)$, where $G$ is defined as
$$G(X,Y,Z):=(AX)*(AY)*(AZ)
  \quad\hbox{for }X,Y,Z\in\CN\,.
  \eqno(5.1)$$
This $G$ is trilinear and symmetric from $\CN\times\CN\times\CN$
into $\CN$, and we can apply Proposition~4.1 to~$F$: for~$\lambda
\in\C\setminus\{0\}$, $|\lambda|\ne1$, there exists a unique
analytic $K_\lambda$, defined as a convergent Taylor series in a
neighbourhood $\dom K_\lambda$ of~$0\in\CN$ and with values in~$\CN$
such that $K_\lambda(0)=0$, $K'_\lambda(0)=I_N$ and such that
$\lambda F(K_\lambda(X))=K_\lambda(\lambda X)$ for all~$X$ such that
$X,\lambda X\in\dom K_\lambda$.

\bigbreak

{\bf Proposition~5.1.} \sl If~$X\in\dom K_\lambda$ and $X_0\in\ker A$
then $X+X_0\in\dom K_\lambda$ and $K_\lambda(X+X_0)=
K_\lambda(X)+X_0$.
\rm

\medbreak

{\bf Proof.} Consider the recursive formulas~(4.2): to start with
$$\Psi_0(X+X_0)=0=\Psi_0(X)\,,\qquad
  \Psi_1(X+X_0)=X+X_0=\Psi_1(X)+X_0\,.
  \eqno(5.2)$$
If $\Psi_r(X+X_0)$ equals either $\Psi_r(X)$ or $\Psi_r(X)+X_0$ for
all $r<m$, then $A\Psi_r(X+X_0)=A\Psi_r(X)$ and
$$\eqalign{
  \Psi_m&(X+X_0)=\cr
  ={}&{1\over1-\lambda^{m-1}}
  \sum_{p+q+r=m \atop 0\le p,q,r<m}
  \bigl(A\Psi_p(X+X_0)\bigr)*
  \bigl(A\Psi_q(X+X_0)\bigr)*
  \bigl(A\Psi_r(X+X_0)\bigr)=\cr
  ={}&{1\over1-\lambda^{m-1}}
  \sum_{p+q+r=m \atop 0\le p,q,r<m}
  \bigl(A\Psi_p(X)\bigr)*
  \bigl(A\Psi_q(X)\bigr)*
  \bigl(A\Psi_r(X)\bigr)=
  \Psi_m(X)
  \quad\hbox{for }m\ge2\,.\cr}
  \eqno(5.3)$$
\qed

\bigbreak

The paired function~$f$ can be written in the form
$f(x)=x-g(x,x,x)$, where $g$ is the trilinear symmetric form
defined by
$$g(x,y,z):=B\bigl((ACx)*(ACy)*(ACz)\bigr)\,.
  \eqno(5.4)$$
Hence Proposition~4.1 can be applied to~$f$ too: for~$\lambda
\in\C\setminus\{0\}$, $|\lambda|\ne1$, there exists a unique
analytic $k_\lambda$, defined as a convergent Taylor series in a
neighbourhood of~$0\in\Cn$ and with values in~$\Cn$ such that
$k_\lambda(0)=0$,
$k'_\lambda(0)=I_n$ and such that $\lambda
f(k_\lambda(x))=k_\lambda(\lambda x)$ for all~$x$ such that
$x,\lambda x\in\dom k_\lambda$.

The next two Propositions teach us that whenever either
$k_\lambda$ or~$K_\lambda$ is globally defined, then the other one
is too, so that the following commutative diagrams always travel
together:
$$
  \matrix{\Cn&\mapleft{k_\lambda}&\Cn\cr
          \mapdown{\lambda f}&&\mapdown{\lambda I_n}\cr
          \Cn&\mapleft{k_\lambda}&\Cn\cr}
  \qquad\qquad
  \matrix{\CN&\mapleft{K_\lambda}&\CN\cr
          \mapdown{\lambda F}&&\mapdown{\lambda I_N}\cr
          \CN&\mapleft{K_\lambda}&\CN\cr}
  \eqno(5.5)$$

\bigskip

{\bf Proposition~5.2.} \sl For small $x\in\Cn$, $X\in\CN$ we have
that $k_\lambda(x)=BK_\lambda(Cx)$ and $k_\lambda(BX)
=BK_\lambda(X)$. Moreover, if $K_\lambda$ is globally defined
on~$\CN$, then the function~$k_\lambda$ is globally defined
on~$\Cn$~too, and the following diagrams commute:
$$
  \matrix{\CN&\mapright B&\Cn\cr
          \mapup{K_\lambda}&&\mapup{k_\lambda}\cr
          \CN&\mapleft C&\Cn\cr}
  \qquad\qquad
  \matrix{\CN&\mapright B&\Cn\cr
          \mapup{K_\lambda}&&\mapup{k_\lambda}\cr
          \CN&\mapright B&\Cn\cr}
  \eqno(5.6)$$
In particular, if~$K_\lambda$ is a polynomial mapping, so
is~$k_\lambda$, and the degree of~$k_\lambda$ is not higher than the
degree of~$K_\lambda$.
\rm

\medbreak

{\bf Proof.} Let $p(x):=BK_\lambda(Cx)$ for small~$x\in\Cn$.
We have that $p(0)=BK_\lambda(0)=0$, $p'(0)=BK'_\lambda(0)C=BC=
I_n$, and
$$\eqalign{
  \lambda f\bigl(p(x)\bigr)={}&
  \lambda BF\bigl(Cp(x)\bigr)=
  \lambda BF\bigl(CBK_\lambda(Cx)\bigr)=\cr
  ={}&
  \lambda B\Bigl(F\bigl(CBK_\lambda(Cx)\bigr)+
  \underbrace{K_\lambda(Cx)-CBK_\lambda(Cx)}_{\hbox to0pt{\hss
    $\scriptstyle\in\ker A=\ker B$\hss}}\Bigr)=\cr
  ={}&
  \lambda BF\Bigl(CBK_\lambda(Cx)+K_\lambda(Cx)-CBK_\lambda(Cx)
  \Bigr)=\cr
  ={}&
  \lambda BF\bigl(K_\lambda(Cx)\bigr)=
  BK_\lambda(\lambda Cx)=\cr
  ={}&p(\lambda x)\,.\cr}
  \eqno(5.7)$$
The function~$p$ is obviously  analytic and it satisfies the same
relations that define~$k_\lambda$ uniquely by Proposition~4.1.
Hence $p=k_\lambda$ near the origin and the conjugation relation
$\lambda f(k_\lambda(x))=k_\lambda(\lambda x)$ holds for small~$x\in
\Cn$. Next, let~$X\in\CN$ be small. From Proposition~5.1 we have that
$$K_\lambda(X)=K_\lambda\Bigl(CBX+
  \underbrace{X-CBX}_{\hbox to0pt{\hss
    $\scriptstyle\in\ker A=\ker B$\hss}}\Bigr)=
  K_\lambda(CBX)+X-CBX\,,
  \eqno(5.8)$$
whence, applying~$B$ we get that $BK_\lambda(X)=BK_\lambda(CBX)=
k_\lambda(BX)$. If~$K_\lambda$ is globally defined, the identities
extend to the whole spaces and define $k_\lambda$ everywhere
on~$\Cn$. The first one shows also that if~$K_\lambda$ is polynomial
so is~$k_\lambda$, with no greater degree.
\qed

\bigbreak

{\bf Proposition~5.3.} \sl For~$X$ in a neighbourhood of~$0\in\CN$ we
can write $K_\lambda(X)= Ck_\lambda(BX)+Q(X)$,
where $Q$ is the unique analytic function such that
$Q'(0)=I_n-CB$ and such that
$$Q(\lambda X)-\lambda Q(X)=
  \lambda(I_n-CB)F\bigl(Ck_\lambda(BX)\bigr)
  \qquad\hbox{for small }X\in\CN\,.
  \eqno(5.9)$$
If $k_\lambda$ is globally defined on~$\Cn$, then $Q$ and~$K_\lambda$
are also globally defined on~$\CN$. Moreover, if~$k_\lambda$
is a polynomial mapping, so is~$K_\lambda$, and the degree
of~$K_\lambda$ is at most three times the degree of~$k_\lambda$.
\rm

\bigbreak

{\bf Proof.} Let $Q$ be defined as $Q(X):=K_\lambda(X)-
Ck_\lambda(BX)$ for small~$X$. This function $Q$ is obviously
analytic near the origin and
$Q'(0)=K'_\lambda(0)-Ck'_\lambda(0)B=I_n-CB$. Using Proposition~5.2
we have that $BQ(X)=BK_\lambda(X)-BCk_\lambda(BX)=k_\lambda(BX)-
k_\lambda(BX)=0$, so that $Q(X)\in\ker B=\ker A$. Let us write the
conjugation relation $\lambda F(K_\lambda(X))=K_\lambda(\lambda X)$
in terms of~$Q$: the left-hand side becomes
$$\lambda F\bigl(Ck_\lambda(BX)+Q(X)\bigr)=
  \lambda F\bigl(Ck_\lambda(BX)\bigr)+\lambda Q(X)\,,
  \eqno(5.10)$$
while the right-hand side is, using the conjugation relation
for~$f,k_\lambda$ and the definition of~$f$,
$$Ck_\lambda(\lambda BX)+Q(\lambda X)=
  \lambda Cf\bigl(k_\lambda(BX)\bigr)+Q(\lambda X)=
  \lambda CBF\bigl(Ck_\lambda(BX)\bigr)+Q(\lambda X)\,.
  \eqno(5.11)$$
Formula~(5.9) is simply the rearranged combination of~(5.10)
and~(5.11). Let $\sum\Phi_m(X)$ be the Taylor expansion of
$X\mapsto (I_n-CB)F(Ck_\lambda(BX))$ centered in the origin (notice
that this function has values in~$\ker A$), and $\sum\varphi_m(X)$
the one of~$Q(X)$. Relation~(5.9) is equivalent to
$$(\lambda^{m-1}-1)\varphi_m(X)=\Phi_m(X)\,,
  \eqno(5.12)$$
which determines uniquely all the terms $\varphi_m$ except the
one with~$m=1$. 

\noindent
If we assume that $k_\lambda$ is globally defined and
$0<|\lambda|<1$, then formula~(5.9) can be used to extend
analytically the definition of~$Q$ from any ball $\{X\;:\;|X|<r\}$ to
the larger ball $\{X\;:\;|\lambda X|<r\}$. This means that $Q$ is
global, and hence~$K_\lambda$ too. When $|\lambda|>1$ both
conjugation are global to begin with, because of Proposition~4.1.

\noindent
If~$k_\lambda$ is a polynomial mapping, then all
the $\Phi_m$ vanish identically for~$m$ beyond three times its
degree, so the same happens for~$\varphi_m$ too.
\qed

\bigbreak

In the remaining part of this Section we will deduce the
invertibility of each of~$K_\lambda,k_\lambda$ from the
invertibility of the other. For this we will assume that
$k_\lambda$ and $K_\lambda$ are both globally defined, as it is
always the case when either $|\lambda|>1$ or $f,F$ are invertible.

\bigskip

{\bf Proposition~5.4.} \sl If $K_\lambda$ is one-to-one, so
$k_\lambda$~is. If~$K_\lambda$ is onto, so $k_\lambda$~is.
If~$K_\lambda$ is bijective, so is~$k_\lambda$, and
$k_\lambda^{-1}(y)=BK_\lambda^{-1}(Cy)$ for all~$y\in\Cn$, i.e.,
the following diagram commutes:
$$
  \matrix{\CN&\mapleft C&\Cn\cr
          \mapdown{K_\lambda^{-1}}&&\mapdown{k_\lambda^{-1}}\cr
          \CN&\mapright B&\Cn\cr}
  \eqno(5.13)$$
In particular, if~$K_\lambda^{-1}$ is a polynomial mapping, then so
is~$k_\lambda^{-1}$.
\rm

\medbreak

{\bf Proof.} Suppose that $K_\lambda$ is one-to-one. Then for
all~$x_1,x_2\in\Cn$
$$\eqalign{
  k_\lambda(x_1)=k_\lambda(x_2)\quad\Longrightarrow{}&\quad
  BK_\lambda(Cx_1)=BK_\lambda(Cx_2)\cr
  \Longrightarrow{}&\quad
  K_\lambda(Cx_1)-K_\lambda(Cx_2)=X_0\in\ker A=\ker B\cr
  \Longrightarrow{}&\quad
  K_\lambda(Cx_1)=K_\lambda(Cx_2+X_0)\cr
  \Longrightarrow{}&\quad
  Cx_1=Cx_2+X_0\cr
  \Longrightarrow{}&\quad
  \range C\ni C(x_1-x_2)=X_0\in\ker A=\ker B\cr
  \Longrightarrow{}&\quad
  C(x_1-x_2)=X_0=0\cr
  \Longrightarrow{}&\quad
  x_1=x_2\,.\cr}
  \eqno(5.14)$$
Suppose that $K_\lambda$ is onto. Let $y\in\Cn$ be
arbitrary. There exists $Y\in\CN$ such that $K_\lambda(Y)=Cy$.
Then
$$\eqalign{
  \range k_\lambda\ni{}&k_\lambda(BY)=
  BK_\lambda\bigl(CBY\bigr)=
  BK_\lambda\Bigl(CBY+
  \underbrace{Y-CBY}_{\hbox to0pt{\hss
    $\scriptstyle\in\ker A$\hss}}\Bigr)=\cr
  ={}&
  BK_\lambda(Y)=BCy=y\,.\cr}
  \eqno(5.15)$$
The inversion formula comes by writing $Y=K_\lambda^{-1}(Cy)$
in~(5.15).
\qed

\bigbreak

{\bf Proposition~5.5.} \sl If~$k_\lambda$ is one-to-one, so
is~$K_\lambda$. If~$k_\lambda$ is onto, so is~$K_\lambda$.
If~$k_\lambda$ is bijective, so is~$K_\lambda$, and
$$K^{-1}_\lambda(Y)=Y-K_\lambda\bigl(Ck_\lambda^{-1}(BY)\bigr)+
  Ck_\lambda^{-1}(BY)
  \quad\hbox{for all }Y\in\CN\,.
  \eqno(5.16)$$
In particular, if~$k_\lambda^{-1}$ is a polynomial mapping, so
is~$K_\lambda^{-1}$, and the degree of~$K_\lambda^{-1}$ is not
larger than the product of the degrees of~$K_\lambda$
and~$k_\lambda^{-1}$.
\rm

\medbreak

{\bf Proof.} Suppose that $k_\lambda$ is one-to-one and let
$X_1,X_2\in\CN$. Then
$$\eqalign{
  K_\lambda(X_1)={}&K_\lambda(X_2)\quad\Longrightarrow\cr
  \Longrightarrow{}&\quad
  K_\lambda\Bigl(CBX_1+
  \underbrace{X_1-CBX_1}_{\hbox to0pt{\hss
  $\scriptstyle\in\ker B=\ker A$\hss}}\Bigr)=
  K_\lambda\Bigl(CBX_2+
  \underbrace{X_2-CBX_2}_{\hbox to0pt{\hss
  $\scriptstyle\in\ker B=\ker A$\hss}}\Bigr)\cr
  \Longrightarrow{}&\quad
  K_\lambda(CBX_1)+X_1-CBX_1=
  K_\lambda(CBX_2)+X_2-CBX_2 \quad\hbox{(*)}\cr
  \Longrightarrow{}&\quad
  BK_\lambda(CBX_1)=BK_\lambda(CBX_2)\cr
  \Longrightarrow{}&\quad
  k_\lambda(BX_1)=k_\lambda(BX_2)\cr
  \Longrightarrow{}&\quad
  BX_1=BX_2\qquad\hbox{(using * above)}\cr
  \Longrightarrow{}&\quad
  X_1=X_2\,.\cr}
  \eqno(5.17)$$
Suppose that~$k_\lambda$ is onto. Let~$Y\in\CN$. There exists
$x\in\Cn$ such that $BY=k_\lambda(x)=BK_\lambda(Cx)$. In
particular $Y-K_\lambda(Cx)\in\ker A=\ker B$. Then
$$\range K_\lambda\ni K_\lambda\Bigl(Cx+
  \underbrace{Y-K_\lambda(Cx)}_{\hbox to0pt{\hss
  $\scriptstyle\in\ker B=\ker A$\hss}}\Bigr)=
  K_\lambda(Cx)+Y-K_\lambda(Cx)=Y\,.
  \eqno(5.18)$$
In particular, if $k_\lambda$ is bijective just write
$x=k_\lambda^{-1} (BY)$ to get the inversion formula. Finally,
if~$k_\lambda^{-1}$ is a polynomial map, then also $k_\lambda$
must be polynomial by a well-known result (see e.g.~[16]), and then
$K_\lambda$ too by Proposition~5.3.
\qed

\bigskip

If we weakened Proposition~5.5 by saying {\sl``if $k_\lambda^{-1}$
and $k_\lambda$ are polynomial mapping, so is~$K_\lambda^{-1}$''},
then we would not need to resort to the advanced complex analysis
result of~[16], and the result would extend to the real case too.

\vfill\eject
\baselineskip=12pt

\centerline{\twelvebf 6. Examples}

\nobreak\bigskip

{\bf Example~6.1.} Consider the $15\times15$ matrix
$$\def\s{\scriptstyle}
  A={1\over2}\left(\matrix{
  \s0&\s0&\s0&\s0&\s0&\s0&\s0&\s0&\s0&\s0 &\s
    0&\s0&\s0&\s0&\s0\cr
  \s0&\s0&\s0&\s0&\s0&\s0&\s0&\s0&\s0&\s0&\s0&
    \s0&\s0&\s0&\s0\cr
  \s0&\s0&\s0&\s-4&\s-2&\s2&\s2&\s2&\s0&\s0&
    \s-2&\s0&\s0&\s-2&\s0\cr
  \s0&\s0&\s-2&\s0&\s-2&\s0&\s1&\s0&\s0&\s1&
    \s0&\s-1&\s-1&\s0&\s0\cr
  \s0&\s0&\s2&\s-4&\s0&\s0&\s0&\s2&\s-2&\s-2&
    \s-2&\s0&\s0&\s0&\s2\cr
  \s2&\s0&\s2&\s-4&\s0&\s0&\s0&\s2&\s-2&\s-2&
    \s-2&\s0&\s0&\s0&\s2\cr
  \s0&\s2&\s2&\s-4&\s0&\s0&\s0&\s2&\s-2&\s-2&
    \s-2&\s0&\s0&\s0&\s2\cr
  \s2&\s0&\s-2&\s0&\s-2&\s0&\s1&\s0&\s0&\s1&
    \s0&\s-1&\s-1&\s0&\s0\cr
  \s2&\s0&\s0&\s-4&\s-2&\s2&\s2&\s2&\s0&\s0&
    \s-2&\s0&\s0&\s-2&\s0\cr
  \s0&\s2&\s0&\s-4&\s-2&\s2&\s2&\s2&\s0&\s0&
    \s-2&\s0&\s0&\s-2&\s0\cr
  \s2&\s0&\s2&\s0&\s2&\s0&\s-1&\s0&\s0&\s-1&
    \s0&\s1&\s1&\s0&\s0\cr
  \s0&\s2&\s-2&\s4&\s0&\s0&\s0&\s-2&\s2&\s2&
    \s2&\s0&\s0&\s0&\s-2\cr
  \s0&\s2&\s0&\s4&\s2&\s-2&\s-2&\s-2&\s0&\s0&
    \s2&\s0&\s0&\s2&\s0\cr
  \s2&\s2&\s2&\s-4&\s0&\s0&\s0&\s2&\s-2&\s-2&
    \s-2&\s0&\s0&\s0&\s2\cr
  \s2&\s2&\s0&\s-4&\s-2&\s2&\s2&\s2&\s0&\s0&
    \s-2&\s0&\s0&\s-2&\s0\cr}\right)\,.
  \eqno(6.1)$$
The function $F\colon\C^{15}\to\C^{15}$ defined by
$F(X)=X-(AX)^{*3}$ was introduced by Dru\.zkowski in~[9] as a
simpler alternative to an example by Rusek~[17], concerning some
geometric condition proposed by Yagzhev.

It can be verified that $A$ has rank equal to~5 and that $A^2=0$. A
linear mapping (or matrix) $B\colon\C^{15}\to\C^5$ with the same
kernel as~$A$ is the following:
$$B={1\over2}\left(
  \matrix{2&0&0&0&0&0&0&0&0&0&0&0&0&0&0\cr
  0&2&0&0&0&0&0&0&0&0&0&0&0&0&0\cr
  0&0&0&-4&-2&2&2&2&0&0&-2&0&0&-2&0\cr
  0&0&2&-4&0&0&0&2&-2&-2&-2&0&0&0&2\cr
  0&0&-2&0&-2&0&1&0&0&1&0&-1&-1&0&0\cr}\right)\,.
  \eqno(6.2)$$
Notice that, if we ignore the first couple of columns, the set
of the rows of~$B$ coincides with the set of the rows of~$A$. It
can be verified that the rows of~$B$ are in fact a basis for the
orthogonal to the kernel of~$A$, with respect to the canonical
scalar product. Anyway, a simple right inverse~$C$ of~$B$ is given
by
$$C^T:=\left(\matrix{
  1&0&0&0&0&0&0&0&0&0&0&0&0&0&0\cr
  0&1&0&0&0&0&0&0&0&0&0&0&0&0&0\cr
  0&0&0&0&0&1&0&0&0&0&0&0&0&0&0\cr
  0&0&0&0&0&0&0&0&-1&0&0&0&0&0&0\cr
  0&0&0&0&0&0&0&0&0&0&0&-2&0&0&0\cr
  }\right)\,.
  \eqno(6.3)$$
The mapping $f\colon\C^5\to\C^5$ paired to~$F$ through~$B$ and~$C$
($f(x):=BF(Cx)$) is calculated as
$$f(x)=x+
  3\left(\matrix{0\cr0\cr
  x_1^2x_2 + x_1x_2^2 + 2x_1x_2x_4 - 2x_1^2x_5\cr
  -x_1^2x_2 - x_1x_2^2 - 2x_1x_2x_3 - 2x_1^2x_5\cr
  -x_2^2x_3 - x_2^2x_4\cr}\right)\,.
  \eqno(6.4)$$
The inverse of~$f$ is easily found by computer and it is a
polynomial mapping of degree~7:
$$\eqalign{
  f^{-1}(y)=y+{}&
  3\left(\matrix{
  0\cr0\cr
  -y_1^2y_2-y_1y_2^2-2y_1y_2y_4+2y_1^2y_5\cr
  \noalign{\smallskip}
  y_1^2y_2 + y_1y_2^2 + 2y_1y_2y_3 + 2y_1^2y_5\cr
  \noalign{\smallskip}
  y_2^2y_3 + y_2^2y_4}\right)+\cr
  \noalign{\smallskip}
  +{}&18\left(\matrix{
  0\cr 0\cr
  -y_1^3y_2^2-y_1^2y_2^3-y_1^2y_2^2y_3+y_1^2y_2^2y_4-
    2y_1^3y_2y_5\cr 
  \noalign{\smallskip}
  -y_1^3y_2^2-y_1^2y_2^3+y_1^2y_2^2y_3-y_1^2y_2^2y_4+
    2y_1^3y_2y_5\cr
  \noalign{\smallskip}
  y_1y_2^3y_3-y_1y_2^3y_4+2y_1^2y_2^2y_5}\right)+\cr
  \noalign{\smallskip}
  +{}&108\left(\matrix{
  0\cr0\cr
  y_1^4y_2^3 + y_1^3y_2^4\cr
  \noalign{\smallskip}
  -y_1^4y_2^3 - y_1^3y_2^4\cr
  \noalign{\smallskip}
  -y_1^3y_2^4 - y_1^2y_2^5}\right)\,.\cr}
  \eqno(6.5)$$
Proposition~3.2 predicts now that the inverse of~$F$ is a
polynomial mapping of degree at most~21 and that it is given by
the formula
$$F^{-1}(Y)=Y+\Bigl(ACf^{-1}(BY)\Bigr)^{*3}=
  2Y-F\bigl(Cf^{-1}(BY)\bigr)\,.
  \eqno(6.6)$$
The pre-conjugation~$k_\lambda$ of the paired mapping~$f$ can be
computed through the recursive formula~(4.2) and turns out to be a
polynomial mapping of degree~7:
$$\eqalign{&k_\lambda(x)=x+
  {3\over1-\lambda^2}
  \left(\matrix{
  0\cr 0\cr
  -x_1^2x_2 - x_1x_2^2 - 2x_1x_2x_4 + 2x_1^2x_5\cr
  x_1^2x_2 + x_1x_2^2 + 2x_1x_2x_3 + 2x_1^2x_5\cr
  x_2^2x_3 + x_2^2x_4\cr}\right)+\cr
  \noalign{\smallskip}
  +{}&
  {18\over(1-\lambda^2)(1-\lambda^4)}
  \left(\matrix{
  0\cr 0\cr
  -x_1^3x_2^2-x_1^2x_2^3-x_1^2x_2^2x_3+x_1^2x_2^2x_4-
    2x_1^3x_2x_5\cr 
  -x_1^3x_2^2- x_1^2x_2^3 + x_1^2x_2^2x_3 - 
    x_1^2x_2^2x_4 + 2x_1^3x_2x_5\cr 
  x_1x_2^3x_3 - x_1x_2^3x_4 + 2x_1^2x_2^2x_5\cr}\right)+\cr
  \noalign{\smallskip}
  +{}&
  {108\over(1-\lambda^2)(1-\lambda^4)(1-\lambda^6)}
  \left(\matrix{
  0\cr 0\cr
  x_1^4x_2^3 + x_1^3x_2^4\cr 
  -x_1^4x_2^3 - x_1^3x_2^4\cr 
  -x_1^3x_2^4 - x_1^2x_2^5\cr}\right)\,.\cr}
  \eqno(6.7)$$
Each homogeneous terms of~$k_\lambda$ is a scalar multiple of the
corresponding term in~$f^{-1}$. This is because the trilinear
form~$g$ associated with~$f$ happens to satisfy
$g(g(x,x,x),g(x,x,x),x)\equiv0$ (see Remark~4.2).

\goodbreak

Using Proposition~5.3 we can predict that the pre-conjugation
$K_\lambda$ for the cubic-linear mapping~$F$ is a polynomial
transformation of degree at most~21. The inverse of~$k_\lambda$
can be computed easily enough, exploiting the fact that $k_\lambda$
is affine in the last three components:
$$\eqalignno{&k^{-1}_\lambda(y)=y+
  {3\over1-\lambda^2}
  \left(\matrix{
  0\cr0\cr
  y_1^2y_2 + y_1y_2^2 + 2y_1y_2y_4 - 2y_1^2y_5\cr
  \noalign{\smallskip}
  -y_1^2y_2 - y_1y_2^2 - 2y_1y_2y_3 - 2y_1^2y_5\cr
  \noalign{\smallskip}
  -y_2^2y_3 - y_2^2y_4\cr}\right)+\cr
  \noalign{\bigbreak}
  +{}&{18\lambda^2\over(1-\lambda^2)(1-\lambda^4)}
  \left(\matrix{
  0\cr0\cr
  -y_1^3y_2^2 - y_1^2y_2^3 - y_1^2y_2^2y_3 + 
     y_1^2y_2^2y_4 - 2y_1^3y_2y_5\cr 
  \noalign{\smallskip}
  -y_1^3y_2^2 - y_1^2y_2^3 + y_1^2y_2^2y_3 - 
    y_1^2y_2^2y_4 + 2y_1^3y_2y_5\cr 
  \noalign{\smallskip}
  y_1y_2^3y_3 - y_1y_2^3y_4 + 2y_1^2y_2^2y_5\cr}\right)+
  &(6.8)\cr
  \noalign{\bigbreak}
  +{}&{108\lambda^6\over(1-\lambda^2)(1-\lambda^4)(1-\lambda^6)}
  \left(\matrix{0\cr 0\cr
  -y_1^4y_2^3 - y_1^3y_2^4\cr 
  \noalign{\smallskip}
  y_1^4y_2^3 + y_1^3y_2^4\cr
  \noalign{\smallskip}
  y_1^3y_2^4 + y_1^2y_2^5\cr}\right)\,.\cr}$$
From Proposition~5.5 we can draw that $K_\lambda$ is invertible
and that $K_\lambda^{-1}$ is a polynomial transformation of degree
at most~$21\cdot7=147$.

\bigbreak\bigskip

{\bf Example~6.2.} The following polynomial mapping of~$\C^4$
$$f(x):=x+\left(\matrix{(x_3x_1+x_4x_2)x_4\cr
  -(x_3x_1+x_4x_2)x_3\cr
  x_4^3\cr
  0\cr}\right)
  \qquad\hbox{for }x=
  \left(\matrix{x_1\cr x_2\cr x_3\cr x_4\cr}\right)\in\C^4\,,
  \eqno(6.9)$$
was introduced by van den Essen in~[10, page~231]. It is a
cubic-homogeneous mapping with polynomial inverse (of degree~7):
$$f^{-1}(y)=y+\left(\matrix{
  -y_1y_3y_4-y_2y_4^2\cr
  y_1y_3^2 + y_2y_3y_4\cr 
  -y_4^3\cr
  0\cr}\right)+
  \left(\matrix{
  y_1y_4^4\cr
  -2y_1y_3y_4^3 - y_2y_4^4\cr
  0\cr0\cr}\right)+
  \left(\matrix{
  0\cr y_1y_4^6\cr0\cr0\cr}\right)\,.
  \eqno(6.10)$$

\goodbreak

It was shown in~[10, page~231] with a very simple degree argument
that the pre-conjugations $k_\lambda$ could not possibly be
themselves polynomial automorphisms. The Taylor series
of~$k_\lambda$ truncated at the degree~7 is
$$\eqalignno{&k_\lambda(x)=x+
  {1\over1-\lambda^2}
  \left(\matrix{
  -x_1x_3x_4- x_2x_4^2\cr
  x_1x_3^2 + x_2x_3x_4\cr 
  -x_4^3\cr
  0\cr}\right)+\cr
  \noalign{\bigbreak}
  +{}&
  {1\over(1-\lambda^2)(1-\lambda^4)}
  \left(\matrix{
  x_1x_4^4\cr
  -2x_1x_3x_4^3 - x_2x_4^4\cr
  0\cr
  0\cr}\right)+&(6.11)\cr
  \noalign{\bigbreak}
  +{}&
  {1\over(1-\lambda^2)(1-\lambda^4)(1-\lambda^6)}
  \Biggl(
  \lambda^2\left(\matrix{
  -x_1x_3x_4^5- x_2x_4^6\cr 
  x_1x_3^2x_4^4 + x_2x_3x_4^5 + x_1x_4^6\cr
  0\cr 0\cr}
  \right)+\left(\matrix{
  0\cr x_1x_4^6\cr 0\cr 0\cr}\right)
  \Biggr)+\cdots\cr}$$
The paper~[10, page~231] left the question open whether $k_\lambda$
was globally defined for~$|\lambda|<1$, and whether it was globally
invertible for~$|\lambda|\ne1$. The problem was later studied in
detail in~[12], and it was found that the pre-conjugations
$k_\lambda$ are in fact {\it analytic automorphisms} of~$\C^4$
for~$|\lambda|\ne1$, and the coefficients of the power series were
also explicitly calculated.

Through the procedure delineated in Proposition~2.1 it is possible
to pair~$f$ to a cubic-linear map $F\colon\C^{16}\to\C^{16}$. The
first step is to write the third-degree part of~$f(x)$ as a sum of
cubes of linear forms, using formulas~(2.1):
$$\def\s{\scriptstyle}
  f(x)-x={1\over24}\left(\matrix{
  \vbox{\halign{$\s#$\hfil\cr
  -8x_2^3+4(x_2 - x_4)^3\cr
  \qquad+(x_1 - x_3 -x_4)^3\cr
    \qquad-(x_1 + x_3 - x_4)^3\cr
    \qquad+4(x_2 + x_4)^3\cr
    \qquad-(x_1 -x_3 +x_4)^3\cr
    \qquad+(x_1 + x_3 + x_4)^3\cr}}\cr
  \noalign{\medskip\hrule\medskip}
  \vbox{\halign{$\s#$\hfil\cr
  8x_1^3-4(x_1- x_3)^3\cr
    \qquad-4(x_1 + x_3)^3\cr
    \qquad-(x_2 - x_3-x_4)^3\cr
    \qquad+(x_2 + x_3 - x_4)^3\cr
    \qquad+(x_2- x_3 +x_4)^3\cr
    \qquad-(x_2 + x_3 + x_4)^3\cr}}\cr
  \noalign{\medskip\hrule\medskip}
  \s x_4^3\cr
  \noalign{\medskip\hrule\medskip}
  \s0\cr}\right)=-
  {1\over24}\left(\matrix{
  \s8&\s0&\s0&\s0\cr
  \s-4&\s0&\s0&\s0\cr
  \s-1&\s0&\s0&\s0\cr
  \s1&\s0&\s0&\s0\cr
  \s-4&\s0&\s0&\s0\cr
  \s1&\s0&\s0&\s0\cr
  \s-1&\s0&\s0&\s0\cr
  \s0&\s-8&\s0&\s0\cr
  \s0&\s4&\s0&\s0\cr
  \s0&\s4&\s0&\s0\cr
  \s0&\s1&\s0&\s0\cr
  \s0&\s-1&\s0&\s0\cr
  \s0&\s-1&\s0&\s0\cr
  \s0&\s1&\s0&\s0\cr
  \s0&\s0&\s-24&\s0\cr}
  \right)^T
  \left(\matrix{
  \s x_2\cr
  \s x_2- x_4\cr 
  \s x_1- x_3- x_4\cr
  \s x_1+ x_3- x_4\cr
  \s x_2+ x_4\cr
  \s x_1- x_3+ x_4\cr 
  \s x_1+ x_3+ x_4\cr
  \s x_1\cr 
  \s x_1- x_3\cr
  \s x_1+ x_3\cr 
  \s x_2- x_3- x_4\cr
  \s x_2+ x_3- x_4\cr
  \s x_2- x_3+ x_4\cr 
  \s x_2+ x_3+ x_4\cr
  \s x_4\cr}
  \right)^{*3}
  \eqno(6.12)$$
Suitable matrices~$B,D,C$ have 16 as the larger size and are given
by
$$\def\s{\scriptstyle}
  \eqalign{
  B={}&
  {1\over24}\left(
  \matrix{
  \s8&\s-4&\s-1&\s1&\s-4&\s1&\s-1&\s0&\s0&\s0&\s0&\s0&\s0&\s0
    &\s0&\s0\cr
  \s0&\s0&\s0&\s0&\s0&\s0&\s0&\s-8&\s4&\s4&\s1&\s-1&\s-1
    &\s1&\s0&\s0\cr
  \s0&\s0&\s0&\s0&\s0&\s0&\s0&\s0&\s0&\s0&\s0&\s0&\s0
    &\s0&\s-24&\s0\cr
  \s0&\s0&\s0&\s0&\s0&\s0&\s0&\s0&\s0&\s0&\s0&\s0
    &\s0&\s0&\s0&\s24\cr}
  \right),\cr
  \noalign{\medbreak}
  D^T={}&\left(
  \matrix{
  \s0&\s0&\s1&\s1&\s0&\s1&\s1&\s1&\s1&\s1&\s0&\s0&\s0&
     \s0&\s0&\s0\cr
  \s1&\s1&\s0&\s0&\s1&\s0&\s0&\s0&\s0&\s0&\s1&\s1&\s1&\s1&\s0
     &\s0\cr
  \s0&\s0&\s-1&\s1&\s0&\s-1&\s1&\s0&\s-1&\s1&\s-1&\s1&\s-1&\s
    1&\s0&\s0\cr
  \s0&\s-1&\s-1&\s-1&\s1&\s1&\s1&\s0&\s0&\s0&\s-1&\s-1
    &\s1&\s1&\s1&\s0\cr}
  \right),\cr
  \noalign{\medbreak}
  C^T={}&\left(\matrix{
  \s3&\s0&\s0&\s0&\s0&\s0&\s0&\s0&\s0&\s0&\s0&\s0&\s0&
    \s0&\s0&\s0\cr
  \s0&\s0&\s0&\s0&\s0&\s0&\s0&\s-3&\s0&\s0&\s0&\s0&\s0&\s0&\s
    0&\s0\cr
  \s0&\s0&\s0&\s0&\s0&\s0&\s0&\s0&\s0&\s0&\s0&\s0&\s0&\s0&\s
    -1&\s0\cr
  \s0&\s0&\s0&\s0&\s0&\s0&\s0&\s0&\s0&\s0&\s0&\s0&\s0&\s0
    &\s0&\s1\cr}\right)
  \cr}
  \eqno(6.13)$$
(the last column of~$B$ and the last row of~$D$ have been added
to make~$B$ of full rank). We will skip writing down a basis
of~$\ker B$ (although it has been used for the computation), and
proceed to the final matrix~$A$:
$$\def\s{\scriptstyle}
  A={1\over24}\left(
  \matrix{
  \s0&\s0&\s0&\s0&\s0&\s0&\s0&\s-8&\s4&\s4&\s1&\s-1&
    \s-1&\s1&\s0&\s0\cr
  \s0&\s0&\s0&\s0&\s0&\s0&\s0&\s-8&\s4&\s4&\s1&\s-1&\s-1&\s1
    &\s0&\s-24\cr
  \s8&\s-4&\s-1&\s1&\s-4&\s1&\s-1&\s0&\s0&\s0&\s0&\s0
    &\s0&\s0&\s24&\s-24\cr
  \s8&\s-4&\s-1&\s1&\s-4&\s1&\s-1&\s0&\s0&\s0
    &\s0&\s0&\s0&\s0&\s-24&\s-24\cr
  \s0&\s0&\s0&\s0&\s0&\s0&\s0&\s-8&\s
    4&\s4&\s1&\s-1&\s-1&\s1&\s0&\s24\cr
  \s8&\s-4&\s-1&\s1&\s-4&\s1&\s-1
    &\s0&\s0&\s0&\s0&\s0&\s0&\s0&\s24&\s24\cr
  \s8&\s-4&\s-1&\s1&\s-4&\s
    1&\s-1&\s0&\s0&\s0&\s0&\s0&\s0&\s0&\s-24&\s24\cr
  \s8&\s-4&\s-1&\s1
    &\s-4&\s1&\s-1&\s0&\s0&\s0&\s0&\s0&\s0&\s0&\s0&\s0\cr
  \s8&\s-4&\s-1
    &\s1&\s-4&\s1&\s-1&\s0&\s0&\s0&\s0&\s0&\s0&\s0&\s24&\s0\cr
  \s8&\s-4&\s-1&\s1&\s-4&\s1&\s-1&\s0&\s0&\s0&\s0&\s0&\s0&
    \s0&\s-24&\s0\cr
  \s0&\s0&\s0&\s0&\s0&\s0&\s0&\s-8&\s4&\s4&\s1&\s-1&\s-1&
    \s1&\s24&\s-24\cr
  \s0&\s0&\s0&\s0&\s0&\s0&\s0&\s-8&\s4&\s4&\s1&\s-1&\s-1&\s1
    &\s-24&\s-24\cr
  \s0&\s0&\s0&\s0&\s0&\s0&\s0&\s-8&\s4&\s4&\s1&\s-1
    &\s-1&\s1&\s24&\s24\cr
  \s0&\s0&\s0&\s0&\s0&\s0&\s0&\s-8&\s4&\s4&\s1
    &\s-1&\s-1&\s1&\s-24&\s24\cr
  \s0&\s0&\s0&\s0&\s0&\s0&\s0&\s0&\s0&\s
    0&\s0&\s0&\s0&\s0&\s0&\s24\cr
  \s0&\s0&\s0&\s0&\s0&\s0&\s0&\s0&\s0
    &\s0&\s0&\s0&\s0&\s0&\s0&\s0\cr}\right).
  \eqno(6.14)$$
It is possible to check that $A^2\ne0$, $A^3=0$.
Through the results of Sections~3 and~5, the cubic-linear
polynomial mapping
$$F(X):=X-(AX)^{*3}
  \qquad\hbox{for }X\in\C^{16}
  \eqno(6.15)$$
is a polynomial automorphism of~$\C^{16}$, and its inverse is of
degree at most~21. The conjugations $K_\lambda$ of~$F$ are {\it
analytic but not polynomial} automorphisms of~$\C^{16}$, for
all~$\lambda\in\C\setminus\{0\}$, $|\lambda|\ne1$.

\vfill\eject

\centerline{\twelvebf 7. References} 

\nobreak\bigskip

\frenchspacing

\item{[1]} Arnol'd V. I., {\it Geometrical methods in the theory of
        ordinary differential equations}, Springer-Verlag, 1983.

\medbreak

\item{[2]} Bass H., Connell E. and Wright D.,
        {\sl The Jacobian conjecture: reduction of
        degree and formal expansion of the inverse},
        Bull. Amer. Math. Soc. {\it 7}, 287--330  (1982).

\medbreak
 
\item{[3]} Bia\l ynicki-Birula A. and Rosenlicht M.,
        {\sl Injective morphisms of real algebraic varieties}, 
        Proc. Amer. Math. Soc. {\it 13}, 200--203  (1962).

\medbreak

\item{[4]} Cima A., van den Essen A., Gasull A.,
        Hubbers E. and Ma\~nosas F.,
        {\sl A polynomial counterexample to the Markus-Yamabe
        conjecture}. Nijmegen Univ., Dept. of Math., Report
        No.~9551 (1995). To appear in Adv. Math.

\medbreak 

\item{[5]} Deng B., Meisters G. H. and Zampieri G.,
        {\sl Conjugation for polynomial mappings},
        Z. angew. Math. Phys. ZAMP {\it 46}, 872--882 (1995).

\medbreak

\item{[6]} Deng B., {\sl Analytic conjugation, global attractor,
        and the Jacobian conjecture}, University of
        Nebra\-ska-Lincoln (1995).

\medbreak

\item{[7]} Dru\.zkowski L. M.,
        {\sl An effective approach to Keller's Jacobian
        conjecture}, Math. Ann. {\it 264}, 303--313  (1983).

\medbreak

\item{[8]} Dru\.zkowski L. M. and Rusek K.,
        {\sl The formal inverse and the Jacobian conjecture},
        Ann. Polon. Math. {\it 46}, 85--90 (1985).

\medbreak

\item{[9]} Dru\.zkowski L. M.,
        {\sl The Jacobian conjecture},
        Institute of Mathematics, Polish Academy of Sciences,
        Preprint {\it 492} (1991).

\medbreak 

\item{[10]} van den Essen  A. (Editor), 
        {\it Automorphisms of Affine Spaces},
        Proceedings of the Cura\c{c}ao Conference, July 4-8, 1994,
        Kluwer Academic Publishers, 1995.

\medbreak

\item{[11]} van den Essen A. and Hubbers E.,
        {\sl Chaotic polynomial automorphisms; counterexamples to
        several conjectures}, 
        Adv. in Appl. Math. {\it18}, 382--388 (1997).

\medskip

\item{[12]} Gorni G. and Zampieri G.,
       {\sl On the existence of global analytic conjugations for
       polynomial mappings of Yagzhev type},
       J. Math. Anal. Appl. {\it201}, 880--896 (1996).

\medskip

\item{[13]} Gorni G. and Zampieri G.,
        {\sl Yagzhev polynomial mappings: on the structure of the
        Taylor expansion of their local inverse},
        Ann. Polon. Math. {\it64}, 285--290 (1996).

\medbreak

\item{[14]} Keller O. H.,
        {\sl Ganze Cremona trasformationen},
        Monatshefte f\"ur Mathematik und Phy\-sik {\it 47},
        299--306 (1939).

\medbreak

\item{[15]} Rabier P. J.,
        {\sl On components of polynomial automorphisms in two
        variables}, Comm. in Algebra {\it 24}, 929--937 (1996).

\medbreak

\item{[16]} Rudin W.,
        {\sl Injective polynomial maps are automorphisms},
        Amer. Math.  Monthly
        {\it 102}, 540--543 (1995).

\medbreak

\item{[17]} Rusek K.,
        {\sl A geometric approach to Keller's Jacobian
        conjecture}, Math. Ann. {\it 264},  315--320 (1983).

\medbreak

\item{[18]} Yagzhev A. V., 
        {\sl Keller's problem},
        Siberian Math. J. {\it 21}, 747--754  (1980).

\bye